\documentclass[11pt]{amsart}
\usepackage{amsmath,amssymb}

\theoremstyle{plain}
\newtheorem{propn}{Proposition}[section]
\newtheorem{thm}[propn]{Theorem}
\newtheorem{lemma}[propn]{Lemma}

\theoremstyle{definition}
\newtheorem{defn}[propn]{Definition}

\theoremstyle{remark}
\newtheorem*{rem}{Remark}

\numberwithin{equation}{section}

\begin{document}

\title[Riemannian foliations of projective space]{Riemannian foliations of projective space admitting complex leaves}

\author{Thomas Murphy}


\date{March 1st, 2013.}

\keywords{Singular Riemannian foliations, complex projective space,  tube formulae, holomorphic curves of complex quadrics.}

\begin{abstract}
Motivated by Gray's work on tubular neighbourhoods of complex submanifolds of the complex projective space $\mathbb{P}^n$ equipped with the Fubini-Study metric, Riemannian
foliations of $\mathbb{P}^n$ are studied.  We  prove that there are no complex Riemannian foliations of
any open subset of $\mathbb{P}^n$ of codimension one. In particular there is no Riemannian foliation of any open subset of the projective plane by Riemann surfaces.  We give an intrinsic proof describing
 how a complex submanifold may arise as an exceptional leaf of a non-trivial singular Riemannian foliation of maximal dimension. Gray's tubular neighbourhood formula is applied to obtain a volume bound for certain holomorphic curves of complex quadrics.
\end{abstract}

\maketitle

\thispagestyle{empty}

\section{Introduction}

\let\thefootnote\relax\footnote{
\emph{Mathematics Subject Classification} Primary  53C40.

\emph{Address}: D\'{e}partment de Math\'{e}matique,
Universit\'{e} Libre de Bruxelles,
 \ Boulevard du Triomphe,
B-1050 Bruxelles,
Belgique.

 \emph{Email}: tmurphy@ulb.ac.be. }

Let $\mathbb{P}^n$ denote the complex projective space  equipped with the standard Fubini-Study metric of holomorphic sectional curvature $4$. This paper studies singular Riemannian foliations of $\mathbb{P}^n$. In particular,  singular Riemannian foliations admitting  a complex leaf have a notably rich and interesting structure.

Singular Riemannian foliations of a Riemannian manifold arise naturally via various constructions.  Associated to each isometric group action by a compact Lie group is a singular Riemannian foliation induced by the orbits of the action on $\mathbb{P}^n$. A second prototypical family of examples may
be constructed from isoparametric functions of $\mathbb{P}^n$. The fundamental question \cite{thor1} is to classify the singular
Riemannian foliations  of a given Riemannian manifold satisfying some natural geometric condition. Two examples: Thorbergsson \cite{thorbergsson} proves that full, irreducible isoparametric foliations of $\mathbb{R}^n$ of codimension greater than two are homogenous, and Lytchak \cite{lytchak} has shown that singular Riemannian foliations do not exist on compact negatively curved manifolds.

We begin with a definition.
\begin{defn}
A foliation $\mathcal{F}$ is a decomposition of a smooth manifold $\overline{M}$ into injectively immersed submanifolds,
called the leaves of the foliation $L$, such that
$$
T_pL = \lbrace X_p : X\in \Xi_{\mathcal{F}} \rbrace,
$$
for every $L\in \mathcal{F}$ and every $p\in L$, where $\Xi_{\mathcal{F}}$ is the space of smooth vector fields that are everywhere tangent to the leaves in $\mathcal{F}$.A  foliation $\mathcal{F}$ on $\overline{M}$ is said to be \emph{Riemannian} if a
geodesic is orthogonal to all or none of the leaves $L$ of $\mathcal{F}$ that it meets.
\end{defn}
The
leaves of maximal dimension are called regular, otherwise they are singular.
$\mathcal{F}$ is  regular if all  leaves have the same dimension, and is otherwise singular. If the foliation is regular one recovers the traditional definition of a foliation.  We will define the real dimension of
a maximal leaf $k$ to be the dimension of the foliation, sometimes writing $\mathcal{F}^k$ to emphasize this. Let
$$S = \bigcup_{i\in \mathcal{I}}S_i$$
denote the disjoint union of singular leaves (i.e. those of non-maximal dimension), where $\mathcal{I}$ is some indexing
set. These will be referred to as the singular leaves.

A submanifold $X \subset \mathbb{P}^n$ is said to be complex if it
 respects the ambient K\"{a}hler structure $J$; that is
$JT_p(X) \subset T_p(X)$ for all $p\in M$.
A singular foliation $\mathcal{F}$ of $\overline{M}$ is said to be \emph{complex} if every leaf $L$ is a complex
submanifold. The  dimension of a complex foliation $\mathcal{F}$ is $dim_{\mathbb{R}}(\mathcal{F})/2$.

For $\mathbb{P}^n$ there are two trivial examples of complex Riemannian foliations; the foliation with only one leaf, and the foliation where each leaf is a point. A less obvious example of a complex Riemannian foliation comes from the fibres of the twistor fibrations $\pi: \mathbb{P}^{2n+1}\rightarrow \mathbb{H}P^n$, where both spaces are equipped with canonical metrics. The leaves of this foliation are totally geodesic $\mathbb{P}^1\subset \mathbb{P}^n$.

Let us state our main results. To motivate our work, observe that it is easy to construct  non-trivial examples of regular complex Riemannian foliations in $\mathbb{C}^n$ of all codimensions. In contrast, we observe the following;\\

\emph{Let $\mathcal{F}$ be a regular Riemannian foliation of $\mathbb{P}^n$. Then $\mathcal{F}$ is complex if, and only if, it is the twistor fibration}.
\\

Proving this involves combining established results of  Escobales \cite{escobales}, Nagy \cite{nagy}, and Wilking \cite{wilking}.
Our first goal is to establish a local version of this fact for the special case where $dim_{\mathbb{C}}(\mathcal{F}) = n-1$.

\begin{thm}\label{thm2} Let $\mathcal{O}\subset \mathbb{P}^n$ be an open subset. Then there is no complex Riemannian foliation $\mathcal{F}$ of $\mathcal{O}$ with $codim_{\mathbb{C}}(\mathcal{F}) = 1$.
\end{thm}
This immediately implies that  any singular foliation of an open set of $\mathbb{P}^2$ by Riemann surfaces, such as those in \cite{brunella}, cannot even locally be Riemannian. It is well-known that there is a large moduli space of such foliations, and they have been the focus of much study.  For example, a holomorphic foliation $\mathcal{F}$ of  $\mathbb{P}^2$ is given by a one form $\omega = \sum_{i=0}^2 A_idz^i$, where $[z_0, z_1, z_2]$ are homogeneous coordinates and the $A_i$  are homogeneous polynomials in $z_0, z_1, z_2$ of degree $k + 1$ such that $\sum_{i=0}^2 z_iA_i$ is identically zero. The singular set of $\mathcal{F}$ is the algebraic variety with equations $A_i = 0$ for $0 \leq i \leq 2$.  The above result is notable as previous rigidity results (such as those of Nagy and Escobales mentioned above) only apply to regular foliations of closed manifolds and so it was unknown whether a holomorphic foliation of $\mathbb{P}^n$ could be Riemannian on some open subset $M\subset \mathbb{P}^n$.

Our next result shows that complex submanifolds of $\mathbb{P}^n$  may only arise as singular leaves of Riemannian foliations of maximal dimension
in very special cases.  Locally one can always construct such a foliation by taking  tubes of sufficiently small radii around any embedded  complex submanifold. Hence the question is when can one extend this to obtain a global foliation of $\mathbb{P}^n$ by taking tubes of larger and larger radius until one ``fills up" the ambient manifold.

\begin{thm}\label{thm3}
A complex submanifold  arises as a singular leaf of a Riemannian foliation $\mathcal{F}^{2n-1}$ on $\mathbb{P}^n$
if, and only if, it is isometric to either
\begin{enumerate}
\item a totally geodesic $\mathbb{P}^k \subset \mathbb{P}^n$ for some $k\in \lbrace 0, \dots, n-1\rbrace$,\\
\item the complex quadric $Q^{n-1} = \lbrace [z]\in \mathbb{P}^n : z_0^2 + \dots + z_n^2 =0 \rbrace \subset
\mathbb{P}^n$,\\
\item the Segre embedding of $\mathbb{P}^1\times \mathbb{P}^k$ into $\mathbb{P}^{2k+1}$,\\
\item the Pl\"{u}cker embedding of the complex Grassmann manifold $G_2(\mathbb{C}^5)$ into $\mathbb{P}^9$, or\\
\item the half spin embedding of $SO(10)/U(5)$ in $\mathbb{P}^{15}$.
\end{enumerate}
\end{thm}

This implies that the foliation is equivalent to those induced by an isometric group action of cohomogeneity one. In
particular, in the case where the singular orbit has maximal dimension this theorem is our answer to a problem posed by Alfred Gray in \cite{graynote}. He asked when the
tubes around a complex hypersurface of $\mathbb{P}^n$ ``fill up" the ambient space. It is not clear what precisely he meant: perhaps  he would have allowed the focal set of the tubes to contain a connected
component of measure zero which is not a submanifold of $\mathbb{P}^n$. It seems most likely to us that he meant when does
a complex hypersurface occur as an exceptional leaf of a singular Riemannian foliation, and this theorem explains how this situation arises.

This result can actually be obtained by using more recent work of Alexandrino \cite{alex}, which was not known when Gray formulated this question. It follows from Alexandrino's work that pulling $\mathcal{F}$ back under the Hopf map gives an isoparametric foliation of $S^{2n+1}$, and therefore $\mathcal{F}$ must itself be isoparametric. Such foliations were classified by Wang \cite{wang}. Our goal is to give an intrinsic direct proof (i.e. not requiring the Hopf map and the work of Alexandrino and Wang) by using Jacobi field theory, along the lines of Gray's approach.

Gray was led to consider this question following a remarkable series of papers \cite{graynote}, \cite{gray2}, \cite{gray3},  where he studied the geometry of complex submanifolds $X$ of
$\mathbb{P}^n$ and obtained a generalization of Weyl's tubular neighbourhood formula for submanifolds of $\mathbb{R}^n$ for such $X\subset \mathbb{P}^n$. Recall that
the minimal focal distance of a submanifold $X\subset \mathbb{P}^n$, $\zeta_{\mathbb{P}^n}(X)$, is defined as the smallest
$r\in \mathbb{R}^+$ such that the tube of radius $r$ around $X$ does not contain a focal point. We note for clarity in what follows the geodesic $C_{\xi}$ is defined by the initial conditions $C_{\xi}(0) = p$, $C_{\xi}' (0) = \xi$.  A focal point is a singularity of the map
$$
exp: \nu(X) \rightarrow \mathbb{P}^n \text{, }  r\xi(p) \rightarrow exp_p(r\xi(p))
$$

where $\nu(X)$ is the normal space of $X$.  We will consider this
number as a crude measure of how ``curved" $X$ is.

Solving the corresponding Riccati equation, this number is intuitvely "inversely related" to the maximum principal curvature (i.e. the maximum eigenvalue of $A_{\xi}(p)$, where  $p\in X$, $\xi \in \nu^1(X)$ is a unit normal vector, and $A_{\xi}(p)$ is the shape operator of $X$. Gray used his  formula to show that the
minimal focal distance of a complex hypersurface $M_d$ may be bounded above by the degree $d$ of the polynomial
cutting out $M_d$. It is striking that this bound is independent of the actual polynomial cutting out $M_d$ and only depends on its degree.

The previous result actually implies that $Q^{n-1}$ and $\mathbb{P}^{n-1}$ are the unique hypersurfaces achieving equality in Gray's bound for the minimal focal distance. Consider the totally geodesic complex quadric $Q^{n-1}$. Let  $\Gamma:\mathbb{P}^1\rightarrow$ $Q^{n-1}$,
$n \geq 3$, be an embedded holomorphic curve. For any normal vector $\xi$ to $Q^{n-1}\subset \mathbb{P}^n$ there is a focal point of $Q^{n-1}$ along $C_{\xi}$ at distance $\pi/4$.

Viewing both objects as a complex submanifolds of $\mathbb{P}^n$, it is natural to try characterize which holomorphic curves of $Q^{n-1}$ are less curved (measured using the minimal focal distance) than  the quadric  in the ambient $\mathbb{P}^n$.  Such curves satisfy the bound  $\zeta_{\mathbb{P}^n}(\Gamma(\mathbb{P}^1) \geq \pi/4$. We conclude by characterizing the holomorphic curves embedded in the quadric satisfying this bound.
\begin{thm} \label{thm4}
If $\zeta_{\mathbb{P}^n}(\Gamma(\mathbb{P}^1))\geq \frac{\pi}{4}$, then
$$
Vol\bigg(\Gamma(\mathbb{P}^1) \bigg) < \frac{\pi.2^{n}}{n-1}.
$$
\end{thm}

Equivalently, $deg(\Gamma(\mathbb{P}^1)) < \frac{2^{n-1}}{n-1}$.

There are holomophic curves of any degree embedded in $Q^{n-1}$ as it is birational to projective space. The point of this theorem is that one would a priori not expect a relationship between the submanifold geometry of a holomorphic curve and its degree. There are analogous results for   embedded holomorphic curves in any complex hypersurface, but in general the formulae become more complicated as the submanifold geometry of the hypersurface does. For holomorphic curves embedded into $\mathbb{P}^{n-1}\subset \mathbb{P}^n$, the same proof shows that the flattest curves are the linear $\mathbb{P}^1\subset \mathbb{P}^{n-1}\subset \mathbb{P}^n$.

\section{Focal set theory}

The  standard approach to focal set theory uses $M$-Jacobi vector fields. We will explain something of this theory, referring the
reader to \cite{bco} for details. Let $M\subset \overline{M}$ be a submanifold, with shape operator $A_{\xi}$, where $\xi$ is a unit normal vector field. For ease of notation in this section we will replace $C_{\xi}$ by $\gamma$. Therefore we let $\gamma : I \rightarrow \overline{M}$ be a geodesic in $\overline{M}$ parameterized by
arc-length, with $0\in I$, $p=\gamma(0)$ and $\dot{\gamma}(0) \in \nu_p(M)$. Suppose $V(s,t) = \gamma_s(t)$ is a smooth
geodesic variation of $\gamma=\gamma_0$ with $c(s)= \gamma_s(0)\in M$ and $\xi(s) = \dot{\gamma}_s(0)\in \nu_{c(s)}M$ for
all $s$. Then the Jacobi field $Y$ along $\gamma$ induced by this geodesic variation may be calculated from the initial
values
$$
Y(0) = \frac{d}{ds}\bigg|_{s=0} V(s,0) = \frac{d}{ds}\bigg|_{s=0}\gamma_s(0) = \frac{d}{ds}\bigg|_{s=0}c(s) = \dot{c}(0)\in
T_pM,
$$
and, from the Weingarten formula,
\begin{align*}
Y'(0) =& \frac{\partial}{\partial s}\bigg|_{s=0}\frac{\partial}{\partial t}\bigg|_{t=0} V(s,t)\\
=& \frac{d}{ds}\bigg|_{s=0}\dot{\gamma}_s(0) = \frac{d}{ds}\bigg|_{s=0} \xi(s) \\
=& -A_{\xi(0)}Y(0)  + \nabla^{\perp}_{Y(0)}\xi.
\end{align*}

Hence $Y$ is a vector field along $\gamma$ satisfying
$$
Y(0) \in T_{\gamma(0)}M \text{ and } Y'(0)+ A_{\dot{\gamma}(0)}Y(0) \in \nu_{\gamma(0)}M.
$$
Such vector fields are called $M$-Jacobi vector fields. As a Jacobi field is determined by its initial values, these form an
$\overline{n}$-dimensional subspace of the $2\overline{n}$-dimensional space of all Jacobi fields along $\gamma$. It is
trivial to see that $t\dot{\gamma}$ is always an $M$-Jacobi vector field, which we disregard as it has no geometric
importance for us.
\begin{defn}
Set $\mathcal{J}(M,\gamma)$ to be the $(\overline{n}-1)$-dimensional vector space of all $M$-Jacobi vector fields along
$\gamma$ perpendicular to the vector field $t\rightarrow t\dot{\gamma}(t)$.
\end{defn}

Let $\nu^1(M)$ be the unit normal sphere bundle. For $r\in \mathbb{R}^+$, set
$$
M_r : = \lbrace exp(r\xi) : \xi\in\nu^1(M) \rbrace.
$$
If $M$ is a closed embedded submanifold, $M_r$ will be a hypersurface for sufficiently small $r$.  It is called the tube of
radius $r$ around $M$. There is a natural splitting
$$
\mathcal{J}(M,\gamma) = \mathcal{J}(M,\gamma)^{\top} \oplus \mathcal{J}(M,\gamma)^{\bot}
$$
of linear subspaces
$$
\mathcal{J}(M,\gamma)^{\top} := \lbrace Y\in\mathcal{J}(M,\gamma): Y'(0) = -A_{\dot{\gamma}(0)}(Y(0)) \rbrace,
$$
and
$$
\mathcal{J}(M,\gamma)^{\perp} := \lbrace Y\in\mathcal{J}(M,\gamma): Y(0) = 0, Y'(0)\in \nu_p(M) \rbrace.
$$
If $M$ is a hypersurface, $ \mathcal{J}(M,\gamma)^{\perp}$ is empty. Denote by $A^r_{\dot{\gamma}(r)}$ the shape operator of
$M_r$.
\begin{defn}
$\gamma(r)$ is a \emph{focal point} of $M$ along $\gamma$ if there exists a non-zero $M$-Jacobi field  $Y \in \mathcal{J}(M,
\gamma)$ with $Y(r)=0$.  The \emph{multiplicity of the focal point} is defined as
$$
dim\lbrace Y\in \mathcal{J}(M,\gamma): Y(r)=0\rbrace.
$$
\end{defn}

$M_r$ is assumed to be a submanifold of $\overline{M}$. Let $\xi$ be a smooth curve in $\nu^1(M)$ with $\dot{\gamma}(0) =
\xi(0)$. Then $V(s,t) = exp(t\xi(s))$ is a smooth geodesic variation of $\gamma$ consisting of geodesics intersecting $M$
perpendicularly. Let $Y$ be the corresponding $M$-Jacobi vector field. Then $Y$ is determined by the initial conditions
$Y(0)=\dot{c}(0)$ and $Y'(0)= \xi(0)$. Here again $c: s \rightarrow V(s,0)\in M$ and $\xi$ is viewed as a vector field
along $c$. Since $\xi$ is of unit length, $Y\in \mathcal{J}(M,\gamma)$. The curve $c_r: s\rightarrow exp(r(\xi(s))$ is
smooth in $M_r$ and hence $Y(r) = \dot{c}_r(0)\in T_{\gamma(r)}M_r$. As any tangent vector of $M_r$ at $\gamma(r)$ arises in
such fashion we have
$$
T_{\gamma(r)}M_r = \lbrace Y(r)| Y\in \mathcal{J}(M,\gamma)\rbrace.
$$
Then one calculates that
\begin{equation}\label{jac}
A^r_{\dot{\gamma}(r)}Y(r) = -(Y'(r))^{\top}.
\end{equation}

This is equivalent along $\gamma$ to the matrix valued \emph{Riccati equation}
$$
A_{\gamma'(r)}' = A_{\gamma'(r)}^2 + K_{\gamma'}(r).
$$
This is the standard method of calculating principal curvatures of tubes around a submanifold. The main technical achievement of this paper is the calculation of the principal curvatures of a principal leaf of a Riemannian foliation where the real codimension of the foliation is greater than one. This is achieved in the special case where the principal leaves are complex hypersurfaces of $\mathbb{P}^n$.

 Let $X\subset \mathbb{P}^n$ be a complex submanifold.  $X$ is obviously minimal, and hence CMC, because it is calibrated.
Alternatively, it is an elementary observation that every complex submanifold of a K\"{a}hler manifold is austere.
Moreover, $X$ is always curvature-adapted. Recall a submanifold $X\subset \mathbb{C}P^n$ is said to be \emph{curvature-adapted} if,
at every point $p\in X$, and $\xi \in \nu_p(X)$,
\begin{enumerate}
\item $K_{\xi}:T_pX \rightarrow T_pX$, and \\
\item $K_{\xi}$ and  $A_{\xi}\oplus Id$ commute.
\end{enumerate}

Here $K_{\xi}(\cdot) = R(\xi, \cdot)\xi$ is the normal Jacobi operator. Such submanifolds have been the focus of
much research since their introduction by d'Atri. We caution here that Gray refers to such submanifolds as compatible in his works. It follows from the definition that a common eigenbasis of $K_{\xi}(p)$ and $A_{\xi}(p)\oplus Id$ exists
$\forall p\in X$, denoted by $U(p)$. Gray's theorem (Theorem 6.14 of \cite{gray0}) then states that the tubes around $X$, $X_r$ are also curvature-adapted, for $r$ sufficiently small. This means the principal curvatures of $X_r$ at the point $C_{\xi}(r)$  are given by solving the following family of ODEs along the geodesic $C_{\xi}(r)$ with respect to the set $U_i(r)$ of common eigenvectors of $K_{\xi}(r)$ and $\tilde{A}_{\xi}(r)$, the shape operator of $X_r$:

$$
 \lambda_i'(r) = \lambda_i^2(r) + \kappa_i^2,
$$
$ i=1, \dots, 2n-1$. Here $\lambda_i(0) = \lambda_i(p)$, $i = 1, \dots, q$, the principal curvature of $U_i(0)\in T_pX$, and  $\lambda_i(0) = -\infty$ where $i = q+1,\dots, 2n$.
Also we have $\kappa_i = 1$ or $2$, because the ambient metric has holomorphic sectional curvature $4$. Gray also shows that the $U_i$ may be assumed to be parallel along $C_{\xi}(r)$. We will assume this throughout.

The curvature adapted hypersurfaces of complex projective space are precisely the \emph{Hopf hypersurfaces}, which are defined as hypersurfaces of $\mathbb{P}^n$ with the property that $A_{\xi}(J\xi) = \alpha(J\xi)$, i.e. the  vector field $J\xi$ tangent to the hypersurface is a principal curvature vector. It is known that $\alpha$ (the Hopf principal curvature) is constant for any Hopf hypersurface.

Whilst one cannot expect to understand the principal curvatures of the shape operator of an arbitrary Riemannian foliation, one can calculate them in the following situation.

Let  $\xi$ denote a unit normal vector at a point $p$ to a regular leaf $L(p) $ of $\mathcal{F}$, a complex Riemannian foliation of codimension one of $\mathcal{O}^{n-1}\subset \mathbb{P}^n$. We assume without loss of generality $C_{\xi}(r) \subset \mathcal{O}$. Extend this to a unit normal vector field of $L(p) \cap \mathcal{O}$. This induces a vector field  $\hat{\xi} \in T_1(\mathcal{O})$ by parallel translation along normal geodesics to $L_0$. It is not hard to see $\hat{\xi}$ is orthogonal to the leaves of $\mathcal{F}$ because  $dim_{\mathbb{C}}(\mathcal{F}) = n-1$.    Along $C_{\xi}(r)$ we also denote $\hat{\xi}(C_{\xi}(r))$ by $\hat{\xi}(r)$. Let $U_i, i=1, \dots, 2n-2$ denote a common eigenbasis of $A_{\hat{\xi}}$ and $K_{\hat{\xi}}$ at $p$. Extend this to a common eigenbasis of both operators on $L(p)\cap \mathcal{O}$. Then extend this by parallel translation along normal geodesics to $\hat{U}$, which at each point of $\mathcal{O}$ is a common eigenbasis of $K_{\hat{\xi}}$ and $A_{\hat{\xi}}$. Again, along $C_{\xi}(r)$ we write $\hat{U}_i(C_{\xi}(r)) = \hat{U}_i(r)$.  In this situation, we have the following:

\begin{propn}
The principal curvatures of the shape operator
$A_{\hat{\xi}}(r)$ along $C_{\xi}(r)$ with respect to $\hat{U}_i$ satisfy  the  following family of ODEs
$$
\lambda_i'(r) =  \lambda_i^2(r) + 1 ,
$$
$ i=1, \dots, 2n-2$ where $\lambda_i(0)$ are the principal curvatures with respect to  $U_i$.
\end{propn}

\proof

Define the tensor $A_{\hat{\xi}} \in T^1_1(\mathcal{O})$ by setting
\begin{enumerate}
\item $A_{\hat{\xi}}(X)(q) = A_{\hat{\xi}}X$, if $X\in \Gamma(TL(q))$\\
\item $A_{\hat{\xi}}(X(q)) = 0$, if $X\in \nu(L(q))$, the normal distribution to the leaf $L(q)$. \\
\end{enumerate}

Denote by $\tilde{A}_{\hat{\xi}}(r)$ the shape operator of $L(p)_r$ at the point $C_{\xi}(r)$. Let $\nabla$ denote the Levi-Civita connection of the Fubini-Study metric on $\mathbb{P}^n$.
We compute that
$$
\langle \nabla_{\hat{U}_i}\hat{\xi}, J\hat{\xi}\rangle(C_{\xi}(r)) = -\langle \hat{U}_i(r), \tilde{A}_{\hat{\xi}}(r)(J\hat{\xi}(r))\rangle = 0
$$
as $L(p)_r$ is a Hopf hypersurface. Similarly

\begin{align*}
\langle \nabla_{[\hat{U}_i,\hat{\xi}]}\hat{\xi}, J\hat{\xi}\rangle (C_{\xi}(r)) &= \langle \nabla_{\nabla_{\hat{U}_i}\hat{\xi}}\hat{\xi}, J\hat{\xi}\rangle(C_{\xi}(r))\\
&= \alpha(r) \langle \nabla_{\hat{U}_i}\hat{\xi}, J\hat{\xi}\rangle(C_{\xi}(r))\\
&= 0
\end{align*}

Therefore, along $C_{\xi}(r)$,

\begin{align*}
\bigg(\nabla_{\hat{\xi}}A_{\hat{\xi}}(\hat{U}_i)\bigg)(C_{\xi}(r)) = & \bigg( \nabla_{\hat{\xi}}(A_{\hat{\xi}}(\hat{U}_i)) - A_{\hat{\xi}}(\nabla_{\hat{\xi}}(\hat{U}_i)) \bigg)(C_{\xi}(r))\\
= & - \bigg( \nabla_{\hat{\xi}}\nabla_{\hat{U}_i}\hat{\xi} -A_{\hat{\xi}}([\hat{\xi}, \hat{U}_i]) + A_{\hat{\xi}}^2(U_i)\bigg)(C_{\xi}(r)\\
= &  \bigg( \hat{U}_i + A_{\hat{\xi}}([\hat{\xi}, \hat{U}_i]) -  A_{\hat{\xi}}([\hat{\xi}, \hat{U}_i])+  A_{\hat{\xi}}^2(\hat{U}_i)\bigg)(C_{\xi}(r))\\
=& \bigg( \hat{U}_i + A_{\hat{\xi}}^2(\hat{U}_i)\bigg)(C_{\xi}(r)) .
\end{align*}

The matrix $K_{\hat{\xi}}$ restricted to $Span\lbrace \hat{U}_i, i=1, \dots, 2n-2\rbrace$  is always the identity with respect to any basis, which yields  the first term of the last line. Finally it is immediate that
$$
(\nabla_{\hat{\xi}} A_{\hat{\xi}})(\hat{U}_i)(C_{\xi}(r)) = \nabla_{\hat{\xi}}(\lambda_i\hat{U}_i)(C_{\xi}(r)) = \nabla_{\hat{\xi}}\lambda_i(C_{\xi}(r)) =: \lambda_i'(r).
$$

Thus the principal curvatures $\lambda_i(r)$ of $\hat{U}_i$ along $C_{\hat{\xi}}(r)$ satisfy
$$
\lambda_i'(r) = \lambda_i^2(r) + 1,
$$
where $\lambda_i(0)$ are the principal curvatures corresponding to $U_i(0)$. \endproof

 \endproof

Setting $\lambda_i(0) = Cot(\theta_i)$, the standard solution to this differential equation is $\lambda_i(r) = Cot(\theta_i - r)$.
Such a Riccati-type equation cannot hold for the principal curvatures of a complex Riemannian foliation of $\mathbb{P}^n$ in general; the twistor fibration is an immediate counterexample. The above argument fails as it is no longer true in general that $\nabla_{\hat{\xi}(r)}(A_{\hat{\xi}}(\hat{U}_i)) = -\nabla_{\hat{\xi}}\nabla_{\hat{U}_i}\hat{\xi}$ if $dim_{\mathbb{C}}(\mathcal{F}) < n-1$.

\section{Proof of Theorem \ref{thm2}.}

To begin, we explain our observation classifying regular Riemannian foliations of projective space with complex leaves. This essentially  follows from  the work of Nagy \cite{nagy}, who proved that a regular Riemannian foliation of an irreducible compact   K\"ahler manifold with complex leaves has either every leaf totally goedesic, or else $\mathcal{F}$ admits  sections. Here we recall that a Riemannian foliation is said to admit sections, or be \emph{polar}, if the horizontal distribution $\mathcal{H}$ of the foliation is integrable. When the ambient space is $\mathbb{P}^n$, then in the first case work of Escobales \cite{escobales} imples $\mathcal{F}$ is one of the twistor fibrations. There is no regular polar foliation of $\mathbb{P}^n$ with complex leaves. If such a foliation were to exist, the section through a regular point $p$ would be a leaf in the \emph{dual foliation} (see \cite{wilking}). Thus the dual foliation would have more than one leaf, a contradiction to a theorem of Wilking \cite{wilking}.   It would be of interest to classify the singular Riemannian foliations of Hermitian symmetric spaces with complex leaves. We will address this question in a forthcoming paper.

\proof Suppose $\mathcal{F}$ is  a complex Riemannian foliation of $\mathcal{O}$, and let $p\in \mathcal{O}$. Denote by $L(q)$ the leaf containing $q$.  At $p\in L(p)^k$, write the non-zero principal curvatures of $L(p)$ with respect to the normal vector field $\hat{\xi}$ in the eigenbasis $E(p)$ as
$Cot(\theta_i)$, $i = 1, \dots, 2k$. Assume without loss of generality $L$ is regular.  Then because $\mathcal{F}$ is
Riemannian one can apply the M-Jacobi theory outlined in the last section to calculate the principal curvatures of $A_{\hat{\xi}}$ at the point $C_{\xi}((r))$ in the terms of the $Cot(\theta_i)$. As $\mathcal{F}$ is a complex Riemannian foliation, we have from the standard solution to the Riccati equation that
$$
\sum_{i=1}^{2k} Cot(\theta_i -t)  = 0
$$
for $|t|<\epsilon$. But then this equation cannot hold as the left-hand side is an increasing function with respect to $t$. The only possibilities remaining are that either there are no normal vectors to the leaf, or the leaves are zero-dimensional.
\endproof

\begin{rem}
Another way to prove this theorem would be to calculate that $\mathcal{F}$ is a polar foliation, which is not hard to see. Then one can use the theory of holonomy Jacobi fields (see \cite{gw}, Chapter 1), which also yields a Riccati equation along $C_{\xi}(r)$ via a similar calculation. One can now compare the principal curvatures in the exact same manner as the proof above. However the proof presented is theoretically far simpler, as one does not need to introduce the concepts of holonomy Jacobi fields and their Riccati equations.
\end{rem}

\section{Proofs of Theorems \ref{thm3} and \ref{thm4}}

Consider the homogeneous equation on $\mathbb{C}^{n}$ given as
$$
\sum_{i=0}^n a_iz_i^d = 0,
$$
where $d$ is a positive integer and $a_i\in \mathbb{C}$. Since it continues to hold when each $z_i$ is replaced by $\lambda
z_i$ for any$\lambda\in \mathbb{C}$, the zero locus of this polynomial is a set of $\mathbb{P}^n$. More generally, let
$P_d(z_0,\dots, z_n)$ be any homogeneous polynomial of degree $d$ and consider the set of zeros of $P_d$:
$$
M_d = \lbrace [z_0,\dots, z_n]\in \mathbb{P}^n : P_d(z_0,\dots, z_n)=0.
$$
Then $M_d$ is a complex submanifold of real codimension two, which we define to be a \emph{complex hypersurface} of degree
$d$.

A remarkable result of Gray \cite{gray2} established the rigidity of embedded complex submanifolds  by establishing an analogue of
Weyl's tubular neighbourhood formula.
Let $X\subset \mathbb{P}^n$ be a complex submanifold of dimension $k$, and let
$\gamma = 1+ \gamma_1 +\dots + \gamma_k$ be the total Chern form. Let $\gamma(t) = 1 + t\gamma_1 + \dots + t^k\gamma_k$ and let $\gamma(t) = \prod_{a=1}^k(1+tx_a)$ be the formal factorization. Denote by $F$ the K\"ahler form of $\mathbb{P}^n$. Denote by $T_X(r)$ the tubular neighbourhood around $X$ of radius $r$.

\begin{thm}\label{g1}(Gray) If $r < \zeta_{\mathbb{P}^n}(X)$, then
$$
Vol\bigg(T_X(r)\bigg) =\frac{1}{n!}\int_X\prod_{a=1}^k\bigg( 1 - \frac{1}{\pi}F + x_a\bigg)\wedge\bigg(\pi.Sin^2(r)+ Cos^2(r)F)^n.
$$
\end{thm}
Implicit in this notation is the convention that all terms which are not of dimension $2k$ are discarded before one integrates. For later use, let us calculate this formula in two special cases.  Let $M_d$ be an embedded complex hypersurface of $\mathbb{P}^n$. If $r < \zeta_{\mathbb{P}^n}(M_d)$,
$$
Vol \bigg(T_{M_d}(r)\bigg) = \frac{\pi^{n}}{(n)!}\bigg(1 - (1 - d.Sin^2(2r))^{n}\bigg).
$$
Let $\Gamma: \mathbb{P}^1 \rightarrow \mathbb{P}^n$ be an embedded holomorphic curve. If $r <
\zeta_{\mathbb{P}^n}\Gamma(\mathbb{P}^1)$, then
\begin{align*}
Vol\bigg(T_{\Gamma(\mathbb{P}^1)}(r)\bigg) = \frac{\pi^{n-1}(Sin^2(r))^{n-1}}{(n-1)!}\bigg(( 1 -
\frac{n+1}{n}Sin^2(r)).&Vol(\Gamma(\mathbb{P}^1)\\ & + \frac{2\pi Sin^2(r)}{n}\bigg).
\end{align*}

Gray used this formula to show that $\zeta_{\mathbb{P}^n}(M_d)$ has an upper bound in terms of the degree $d$  of the polynomial cutting out $M_d$. Specifically, he showed
$$
\zeta_{\mathbb{P}^n}(M_d) \leq Sin^{-1}\bigg(\frac{1}{\sqrt{d}}\bigg).
$$
Motivated by this, he then asked when the tubes around a complex hypersurface fill up $\mathbb{P}^n$.

Before we begin the proof, we have to define some notation. For $X\subset \mathbb{P}^n$ complex, we know that the tubes of sufficiently small radius $X_r$ are again curvature-adapted. In the situation we are interested in, the tubes actually foliate projective space. There is a second singular leaf as $\mathbb{P}^n$ is positively curved, and it is shown in \cite{bv} that this leaf is also curvature-adapted. Therefore one can also calculate the principal curvatures in the tube around this second singular leaf using the same techniques as Lemma 7.8 of \cite{gray0}. Travelling along the geodesic $C_{\xi}$ we again reach $X$ and so this gives us two ways of calculating the principal curvatures of $X$; this is the crucial idea in the proof.

The Riccati equation for hypersurfaces is equivalent to the Jacobi equation, where there are also  well-established techniques to describe the principal curvatures of nearby parallel hypersurfaces. Recall that a focal point of a curvature-adapted hypersurface $M$ along $C_{\xi}$ is given as $C_{\xi}(r_0)$ when $J(r_0)=0$ for an $M$-Jacobi vector field. If $M$ is a curvature-adapted hypersurface, this is equivalent to a principal curvature function developing a singularity (i.e. for one of the principal curvature functions $\lambda_i(r)$ corresponding to an eigenvector in the common eigenbasis $U(C_{\xi}(r))$, one has $\lambda_i(r_0) = \infty$). In this case, we say the corresonding principal curvature blows up. The focal set is the union of focal points of a curvature-adapted hypersurface. In our situation,  this is a disjoint set of two embedded submanifolds, one of which is the complex submanifold $X\subset \mathbb{P}^n$.

\begin{lemma}
\label{tm2} Let $\mathcal{F}^{2n-1}$ be a Riemannian foliation on $\mathbb{P}^n$. If all singular leaves $S_i$ are complex, they are locally isometric to totally geodesic $\mathbb{P}^k, k <n$.
\end{lemma}
\proof  As $\mathcal{F}$ is Riemannian, the distance between two singular leaves,
$S_1$ and $S_2$, along $C_{\xi}$ is a constant $r$. Since $S_1$ is complex, the  principal curvature corresponding to $J\xi$
is calculated using Jacobi field theory to be $2Cot(2r)$, where $r$ denotes the distance from $S_1$. Thus $r =\frac{\pi}{2}$, as otherwise $S_2$ would not be complex, a contradiction. Suppose $S_1$ is not totally geodesic. If $S_1$ is
not totally geodesic then there exists a $\delta = Cot(\theta) >0$ such that $\delta, -\delta$ are principal curvatures at
some point $p\in S_1$, for $0<\theta< \frac{\pi}{2}$. The corresponding principal curvature functions for the hypersurface
at distance $r$ are seen to be $Cot(\theta - r)$ and $Cot(-\theta - r)$. But if $0<\delta < \frac{\pi}{2}$ then $Cot(\theta
- r)$ blows up before $t = \frac{\pi}{2}$, a contradiction.
\endproof

\begin{defn}
A connected submanifold  of $\mathbb{P}^n$ is said to have constant principal curvatures if, for every  unit normal vector field
$\xi$ the eigenvalues of the shape operator
$A_{\xi}$ are constant at every point.
\end{defn}
In particular the principal curvatures are independent of the point and normal direction.  This is perhaps not the most natural definition: more usual is the definition used by Heintze, Olmos and Thorbergsson \cite{hot}, but it is most useful for our purposes. Kimura \cite{kimura} gave an explicit classification of the complex submanifolds with constant principal curvatures, which is precisely the list in our Theorem \ref{thm3}.

We now prove Theorem \ref{thm3}.

\proof  $\mathcal{F}$ is induced by tubes around $X$.   We already noted that $S$
consists of two connected components, $X=S_1$ and $S_2$. Picking a point $p\in X$ and looking at a section $C_{\xi}(r)$
passing through $p$ and a corresponding point $q\in S_2$, it follows that $\frac{\pi}{g}, g\in \mathbb{N}$ is the distance
between the two focal sets. It may be assumed that $g > 2$ by Lemma \ref{tm2}. Then $\alpha(r)$, the principal curvature function corresponding to $J\xi(r)$, is given as $2Cot(2r)$. Setting $t=\frac{2\pi}{g}$ one obtains
$$
-\infty = Cot(4\frac{\pi}{g})
$$
and so $g=4$.  Solving the Riccati equation corresponding to the  remaining principal curvatures functions whose principal
curvature vectors span $T_pX$ gives $Cot(\theta_i - r)$.

If $X$ is not totally geodesic the positive principal
curvatures $\lambda_i = Cot(\theta_i)$, $i=2, \dots, \frac{k}{2}, k < n$ have $\frac{\pi}{4} \leq \theta_i \leq \frac{\pi}{2}$. If $\theta_i> \frac{\pi}{4}$, then $\lambda_i$ does not focalize at $S_2$ unless $\theta_i = \frac{\pi}{2}$.
Otherwise it must focalize at $C_{\xi}(\frac{\pi}{4})$, a contradiction.  But it is impossible for $\theta_i <\frac{\pi}{4}$ to hold as the corresponding
principal curvature function would then focalize before $r=\frac{\pi}{4}$. Hence all non-zero principal curvatures of $X$
are $\pm 1$ or $0$. Hence $X$ has constant principal curvatures.
\endproof

Finally we establish Theorem \ref{thm4}.

\proof Consider  the tubular neighbourhood of radius $r$ around $\Gamma(\mathbb{P}^1)\subset \mathbb{P}^n$ , $T_{\Gamma(\mathbb{P}^1)}(r)$
and the tubular neighbhourhood of radius $r$ around $Q^{n-1}\subset \mathbb{P}^n$, $T_{Q}(r)$. If $r < \pi/4$ then
$$
Vol\bigg(T_{\Gamma(\mathbb{P}^1)}(r)\bigg) < Vol\bigg(T_{Q}(r)\bigg).
$$
Applying Theorem \ref{g1} and simplifying, this becomes
$$
0 < Vol\bigg(\Gamma(\mathbb{P}^1)\bigg) < \bigg(\frac{\pi}{n}\bigg)\frac{(1- (1- 2Sin^2(r))^n}{(Sin^2(r))^{n-1} -
\frac{n+1}{n}(Sin^2(r))^n}.
$$
In the limit as $r \rightarrow \frac{\pi}{4}$ the result follows.

\endproof

Obviously, one can also show an analogous bound for holomorphic curves passing through the point $p\in M_d$ of any smooth hypersurface corresponding to the maximum eigenvalue of the shape operator $A_{\xi}, \xi \in \nu^1(M_d)$.

\section{acknowledgments}
This research was supported by an ERC Grant. I wish to acknowledge the helpful advice of  J\"urgen Berndt and the hospitality of the Department of Mathematics at King's College London, England, where part of this paper was written as a Visiting Researcher. I also wish to thank Gudlaugar Thorbergsson for alerting me to Wilking's work, and also Steve Hurder and Alexander Lytchak for helpful comments.

\end{document}